\documentclass[a4paper,12pt]{article}
\usepackage{amsmath,amssymb,amsthm}
\usepackage[all]{xy}
\newtheorem{dfn}{Definition}[section]
\newtheorem{prop}[dfn]{Proposition}
\newtheorem{thm}[dfn]{Theorem}
\newtheorem{lem}[dfn]{Lemma}

\begin{document}
\title{A combinatorial formula for Earle's twisted 1-cocycle on
the mapping class group $\mathcal{M}_{g,*}$}
\author{Yusuke Kuno}
\date{}
\maketitle

\begin{abstract}
We present a formula expressing Earle's twisted 1-cocycle
on the mapping class group of a closed oriented surface of genus $\ge 2$
relative to a fixed base point, with coefficients in the first homology group of the surface.
For this purpose we compare it with Morita's twisted 1-cocycle which is combinatorial.
The key is the computation of these cocycles on a particular element of the
mapping class group, which is topologically a hyperelliptic involution.
\end{abstract}

\noindent \textbf{Introduction and statement of the result.}
Let $(\Sigma_g,*)$ be a closed oriented $C^{\infty}$-surface of genus $\ge 2$
with a fixed base point $*$ and let $\mathcal{M}_{g,*}$ be the mapping class group of
$(\Sigma_g,*)$, namely the group of all orientation preserving diffeomorphisms
of $(\Sigma_g,*)$ modulo isotopies fixing the base point $*$.
The group $\mathcal{M}_{g,*}$ naturally acts on the first homology group $H=H_1(\Sigma_g;\mathbb{Z})$.

In \cite{E}, C. Earle discovered a twisted 1-cocycle
$\psi\colon \mathcal{M}_{g,*}\rightarrow \frac{1}{2g-2}H$.
This cocycle is complex analytic by nature. In fact, he discovered this cocycle in the study
of the action of $\mathcal{M}_{g,*}={\rm mod}(\Gamma)$ on $J(V)$, using his notation, the family of
Jacobi varieties over the Teichm\"uller space of compact Riemann surfaces of genus $g$.
We call $\psi$ \textit{Earle's twisted 1-cocycle}. The construction of $\psi$ will be recalled in section 2.

In view of \cite{M}, $\psi$ gives rise to a generator of the first cohomology group
$H^1(\mathcal{M}_{g,*};H)\cong H^1(\mathcal{M}_{g,1};H)\cong \mathbb{Z}$.
Here $\mathcal{M}_{g,1}$ is the mapping class group of $\Sigma_g$ relative to an
embedded disc, see section 1. Amusingly, other than $\psi$ there have been known various ways of constructing
cocycles representing a generator of this cohomology group; see S. Morita \cite{M,M2,M3,M4}
and T. Trapp \cite{T}. Among others, there is a combinatorial one:
\textit{Morita's twisted 1-cocycle} $f\colon \mathcal{M}_{g,*}\rightarrow H$ defined in \cite{M}. 

Although $\psi$ naturally arises it seems more abstract than other known cocycles.
At first glance, the construction of $\psi$ does not tell much about its value on
a given element of $\mathcal{M}_{g,*}$ which is, for example, expressed as a product
of Dehn twists. The aim of the present paper is to improve this unsatisfactory situation.
For this purpose we compare $\psi$ with $f$.
As a product we obtain a formula expressing Earle's cocycle $\psi$, which appeared
in the context of complex analysis, using Morita's cocycle $f$, more combinatorial one.

To state the result, let us fix the notation.
Let $A_1,\ldots,A_g,B_1,\ldots,B_g$ be a fundamental system of generators
of the fundamental group $\pi_1(\Sigma_g,*)$ and fix it throughout this paper.
Then the group $\pi_1(\Sigma_g,*)$ is isomorphic to the group
$$\Gamma=<A_1,\ldots,A_g,B_1,\ldots,B_g | \zeta=1>.$$
Here, $\zeta=\prod_{k=1}^g [A_k,B_k]=[A_1,B_1]\cdots[A_g,B_g]$ and $[A_k,B_k]=A_kB_kA_k^{-1}B_k^{-1}$.
The natural projection
$\theta\colon \Gamma \rightarrow H$ gives the abelianization of $\Gamma$ and
$H$ is identified with $\mathbb{Z}^{2g}$ by the direct decomposition
$$H=\mathbb{Z}\cdot \theta(A_1)\oplus \cdots \oplus \mathbb{Z}\cdot \theta(B_g).$$

We denote the action of $\mathcal{M}_{g,*}$ on $H$ by $\rho\colon \mathcal{M}_{g,*}\rightarrow Sp(H)$.
Here $Sp(H)$ denotes the group of automorphisms of $H$ preserving the intersection form.
In view of the identification $H\cong \mathbb{Z}^{2g}$ given above, $Sp(H)$ is identified with
the symplectic group $Sp(2g;\mathbb{Z})$.
Let $a_0\in \frac{1}{2g-2}\mathbb{Z}^{2g}$ be the column vector defined by
$$a_0=\frac{1}{2g-2}\ ^{t}(\underbrace{0,\ldots,0}_{g},\underbrace{2,\ldots,2}_{g}).$$
The result of this paper is:

\begin{thm}
\label{main.thm}
We have
$$\psi=-\frac{1}{2g-2}f+\delta a_0,$$
where the coboundary $\delta a_0$ is given by $\delta a_0(\phi)=\rho(\phi)^{-1}\cdot a_0-a_0,
\ \phi\in \mathcal{M}_{g,*}$.
\end{thm}

The organization of this paper is as follows. In section 1 we introduce some groups other
than $\mathcal{M}_{g,*}$ and recall a result on $H^1(\mathcal{M}_{g,*};H)$ by Morita.
Sections 2 and 3 are devoted to the review of Earle's cocycle and Morita's cocycle, respectively.
In section 4 a particular element $\bar{\iota}\in \mathcal{M}_{g,*}$, which is topologically a
hyperelliptic involution, is defined.
Theorem \ref{main.thm} will be proved in section 5.
The key lemma to prove Theorem \ref{main.thm} is Lemma \ref{lem:5-2} in which we compute
the value $\psi(\bar{\iota})$ explicitly.
This computation will be performed in section 6, in which we also prove Lemma \ref{lem:5-1},
another lemma needed to the proof of Theorem \ref{main.thm}.

\section{Mapping class groups, $H^1(\mathcal{M}_{g,*};H)$}
We first introduce some groups related to $\mathcal{M}_{g,*}$.
Let $D\subset \Sigma_g$ be an embedded closed 2-disk centered at $*$ and
let $\mathcal{M}_{g,1}$ be the mapping class group of $(\Sigma_g,D)$,
namely the group of all orientation preserving diffeomorphisms
of $(\Sigma_g,D)$ modulo isotopies fixing $D$.
Choose a base point $*^{\prime}$ of $\Sigma_g \setminus {\rm Int}D$ in $\partial D$.
The fundamental group $\pi_1(\Sigma_g \setminus {\rm Int}D,*^{\prime})$ is
isomorphic to a free group of rank $2g$.
By an appropriate choice of (homotopy class of) based loops
$A_1^{\prime},\ldots,A_g^{\prime},B_1^{\prime},\ldots,B_g^{\prime}$
in $\Sigma_g \setminus {\rm Int}D$, we can arrange that;
\begin{enumerate}
\item $\pi_1(\Sigma_g \setminus {\rm Int}D,*^{\prime})$ is freely generated by
$A_1^{\prime},\ldots,A_g^{\prime},B_1^{\prime},\ldots,B_g^{\prime}$.
\item $\zeta^{\prime}=\prod_{k=1}^g [A_k^{\prime},B_k^{\prime}]$ is homotopic ${\rm rel}\ *^{\prime}$
to the boundary loop $\partial D$.
\item Let $\gamma$ be a path in $D$ joining $*$ and $*^{\prime}$.
Joining with $\gamma$, we have a natural homomorphism from $\pi_1(\Sigma_g \setminus {\rm Int}D,*^{\prime})$
to $\pi_1(\Sigma_g,*)$. Then for each $k$, $A_k^{\prime}$ (resp. $B_k^{\prime}$) is mapped to $A_k$ (resp. $B_k$). 
\end{enumerate}
Henceforth we write $A_k$ (resp. $B_k$) instead of $A_k^{\prime}$ (resp. $B_k^{\prime}$) for
simplicity. Let $F=F(A_1,\ldots,A_g,B_1,\ldots,B_g)$ be the free group generated by $A_1,\ldots,A_g,B_1,\ldots,B_g$.
The mapping class groups $\mathcal{M}_{g,1}$ and $\mathcal{M}_{g,*}$ act on the fundamental groups of the surfaces.
By the well-known theorem of Dehn-Nielsen, these actions are faithful and we have the isomorphisms
$$\mathcal{M}_{g,1}\stackrel{\cong}{\rightarrow} \{ \phi \in Aut(F); \phi(\zeta)=\zeta \}$$
and
$$\mathcal{M}_{g,*}\stackrel{\cong}{\rightarrow} Aut^+(\Gamma),$$
where $\ ^+$ means acting
on the second homology $H_2(\Gamma)\cong \mathbb{Z}$ as the identity.
We identify $\mathcal{M}_{g,1}$ (resp. $\mathcal{M}_{g,*}$) with its image
of the above isomorphism. Note that these identifications depend on the choice of
a fundamental system of generators of $\pi_1(\Sigma_g,*)$.
We also consider the subgroup $\mathcal{N}$ of $Aut(F)$, including $\mathcal{M}_{g,1}$, defined by
$$\mathcal{N}:=\{ \phi\in Aut(F); \phi(\zeta)\ {\rm is \ conjugate \ to}\ \zeta \}.$$
Each element of $\mathcal{N}$ induces an automorphism of $\Gamma$ acting on
$H_2(\Gamma)$ as the identity. Thus we have a homomorphism
$$\pi \colon \mathcal{N}\rightarrow Aut^+(\Gamma)\cong \mathcal{M}_{g,*},$$
and it is well known that the restriction of $\pi$ to $\mathcal{M}_{g,1}$ gives rise to
the central extension
\begin{equation}
\label{eq:1-1}
0 \rightarrow \mathbb{Z} \rightarrow \mathcal{M}_{g,1} \stackrel{\pi|_{\mathcal{M}_{g,1}}}{\rightarrow}
\mathcal{M}_{g,*}\rightarrow 1.
\end{equation}
Topologically, $\pi|_{\mathcal{M}_{g,1}}$ is induced by regarding diffeomorphisms of
$(\Sigma_g,D)$ as diffeomorphisms of $(\Sigma_g,*)$, and the generator of the kernel corresponds
to the Dehn twist along the boundary $\zeta=\partial D$.

We next recall a result on $H^1(\mathcal{M}_{g,*};H)$. First, let us introduce the conventions
in this paper. Let $G$ be a group and $M$ a (left) $G$-module.
By a \textit{twisted 1-cocycle} is meant a map $\Phi\colon G\rightarrow M$ satisfying
$$\Phi(\gamma_1\gamma_2)=\gamma_2^{-1}\cdot \Phi(\gamma_1)+\Phi(\gamma_2)$$
for all $\gamma_1,\gamma_2 \in G$. For $m\in M$, by the \textit{coboundary of $m$} is meant the
map $\delta m\colon G\rightarrow M$ defined by
$$\delta m(\gamma)=\gamma^{-1}\cdot m-m.$$
As usual coboundaries are twisted 1-cocycles. The quotient
$$\{ {\rm twisted} \ 1{\rm -cocycles} \}/ \{ {\rm coboundaries} \}$$
is denoted by $H^1(G;M)$
and called the \textit{first cohomology group of $G$ with coefficients in $M$}.

Regard the group $\pi_1(\Sigma_g,*)$ as a subgroup of $\mathcal{M}_{g,*}$
via inner automorphism; For $x\in \pi_1(\Sigma,*)\cong \Gamma$,
the map $\Gamma \rightarrow \Gamma, y\mapsto xyx^{-1}$ is the corresponding element of
$\mathcal{M}_{g,*} \cong Aut^+(\Gamma)$. In \cite{M}, Morita determined the cohomology group
$H^1(\mathcal{M}_{g,*};H)$:

\begin{prop}[Morita \cite{M}]
\label{prop:1-1}
The cohomology group $H^1(\mathcal{M}_{g,*};H)$ is isomorphic to the infinite cyclic
group $\mathbb{Z}$. A twisted 1-cocycle $\Phi \colon \mathcal{M}_{g,*}\rightarrow H$
represents a generator of $H^1(\mathcal{M}_{g,*};H)$ if and only if the restriction
of $\Phi$ to $\pi_1(\Sigma_g,*)$ coincides with $\pm (2g-2)$ times the abelianization:
$\Phi|_{\pi_1(\Sigma_g,*)}=\pm (2g-2)\theta$.
\end{prop}

\section{Earle's twisted 1-cocycle}
We review Earle's twisted 1-cocycle $\psi$ described in \cite{E}.
Let $\mathcal{T}_{g,1}$ be the Teichm\"uller space of
compact Riemann surfaces of genus $g$ with one distinguished point;
It is the set of equivalence classes of all triads $(X,p,f)$ such that
\begin{enumerate}
\item $X$ is a compact Riemann surface of genus $g$,
\item $f\colon (\Sigma_g,*)\rightarrow (X,p)$ is an orientation preserving diffeomorphism,
\end{enumerate}
and the equivalence relation is defined as follows:
$(X,p,f)$ and $(X^{\prime},p^{\prime},f^{\prime})$ are equivalent if
there exists a biholomorphic map $h\colon (X,p)\rightarrow (X^{\prime},p^{\prime})$ such that
$h\circ f$ is homotopic to $f^{\prime}$.
We denote by $[X,p,f]$ the point of $\mathcal{T}_{g,1}$ represented by $(X,p,f)$.
The mapping class group $\mathcal{M}_{g,*}$ acts on $\mathcal{T}_{g,1}$ by
$$\phi \cdot [X,p,f]=[X,p,f\circ \phi^{-1}],$$
where $\phi \in \mathcal{M}_{g,*}$ and $[X,p,f]\in \mathcal{T}_{g,1}$.
$\mathcal{T}_{g,1}$ has the natural complex structure and is isomorphic to the Bers fiber space.

Now the holomorphic map $\eta\colon \mathcal{T}_{g,1} \rightarrow \mathbb{C}^g$ is defined as follows.
Let $[X,p,f]$ be an element of $\mathcal{T}_{g,1}$. Then $f_*A_1,\ldots,f_*A_g,f_*B_1,\ldots,f_*B_g$ is a
fundamental system of generators of $\pi_1(X,p)$. Let $\omega_1,\ldots,\omega_g$ be the basis of holomorphic
1-forms on $X$ satisfying the normalized condition
$$\int_{[f_*A_i]}\omega_j=\delta_{ij},$$
where $\delta_{ij}$ is the Kronecker delta and $[T]\in H_1(X;\mathbb{Z})$ denotes the homology class 
represented by the loop $T$.
Set $\tau_{ij}=\int_{[f_*B_i]}\omega_j$. The $g\times g$ matrix $\tau=(\tau_{ij})_{i,j}$ is the period matrix of
$X$ with respect to the symplectic basis $[f_*A_1],\ldots,[f_*A_g],[f_*B_1],\ldots,[f_*B_g]$.

For each pair $(j,k)$ $(1\le j,k\le g)$, the quadratic period class corresponding to $\omega_j$ and $\omega_k$
is the function $Q_{jk}\colon \pi_1(X,p)\rightarrow \mathbb{C}$ defined by
$$Q_{jk}(T):=\int_T\omega_k\omega_j=\int_{s=0}^1 \omega_k(T(s))\int_{u=0}^s \omega_j(T(u)),$$
where $T\in \pi_1(X,p)$. Finally, define $\eta=\ ^t(\eta_1,\ldots,\eta_g)$ by
$$(1-g)\eta_j([X,p,f]):=-\frac{1}{2}\tau_{jj}+\sum_{k=1}^g Q_{jk}(f_*A_k),\ 1\le j\le g.$$
It is known that $(1-g)\eta([X,p,f])$ is the vector of Riemann constants for $(X,p)$ with respect to the symplectic
basis $[f_*A_1],\ldots,[f_*A_g],[f_*B_1],\ldots,[f_*B_g]$.

For $\phi\in \mathcal{M}_{g,*}$, we write $\rho(\phi)^{-1}=\left( \begin{array}{cc}
a & b \\
c & d \\ \end{array} \right)$ where $a,b,c$, and $d$ are $g\times g$ integral matrices. 
Set
$$A=A(\phi,[X,p,f]):=(a+\tau c)^{-1}\in GL(g,\mathbb{C})$$ and
\begin{equation}
\label{eq:2-1}
w=w(\phi,[X,p,f]):=A^{-1}\cdot \eta(\phi \cdot [X,p,f])-\eta([X,p,f])\in \mathbb{C}^g.
\end{equation}

Then Theorem 6.6 in \cite{E} says that the vector $w$ can be decomposed into
the factor from the Teichm\"uller space part and the factor from the mapping class group part;
Namely there exists the uniquely determined vector $\psi(\phi)\in \frac{1}{2g-2}\mathbb{Z}^{2g}$
satisfying $w=(I,\tau)\psi(\phi)$. Here $I$ is the $g\times g$ identity matrix.
In this way we obtain a map $\psi\colon \mathcal{M}_{g,*}\rightarrow \frac{1}{2g-2}\mathbb{Z}^{2g}$.

Moreover, Earle showed the following:

\begin{prop}[Earle \cite{E}]The map $\psi$ is a twisted 1-cocycle:
$$\psi(\phi_1\phi_2)=\rho(\phi_2)^{-1}\psi(\phi_1)+\psi(\phi_2),\ \phi_1,\phi_2 \in \mathcal{M}_{g,*}.$$
Further, when restricted to the subgroup $\pi_1(\Sigma_g,*)\subset \mathcal{M}_{g,*}$, $\psi$
coincides with the abelianization: $\psi|_{\pi_1(\Sigma_g,*)}=\theta$.
\label{prop:2-1}
\end{prop}

In this paper we call
$\psi\colon \mathcal{M}_{g,*} \rightarrow \frac{1}{2g-2}\mathbb{Z}^{2g}\cong \frac{1}{2g-2}H$
\textit{Earle's twisted 1-cocycle}. 

\section{Morita's twisted 1-cocycle}
We review Morita's twisted 1-cocycle defined in \cite{M}, section 6.
The abelianization $F^{ab}$ can be naturally identified with $H$ hence inherits the intersection
form from $H$. For $x\in F$, We denote by $[x]$ the element of $F^{ab}=H$ represented by $x$. 

Let $F(\alpha,\beta)$ be the free group generated by $\alpha$ and $\beta$.
For $i=1,\ldots,g$, let $p_i$ be the homomorphism from $F$ to $F(\alpha,\beta)$
defined by $p_i(A_j)=p_i(B_j)=1$ if $j\neq i$, and $p_i(A_i)=\alpha$, $p_i(B_i)=\beta$. 
Any element $x\in F(\alpha,\beta)$ can be uniquely written as the form
$$x=\alpha^{\varepsilon_1}\beta^{\delta_1}\cdots \alpha^{\varepsilon_n}\beta^{\delta_n},$$
where $\varepsilon_i,\delta_i \in \{ -1,0,1 \}$.
We first set
\begin{equation}
\label{eq:3-1}
d(x):=\sum_{k=1}^n\varepsilon_k \sum_{\ell=k}^n\delta_{\ell}
-\sum_{k=1}^n\delta_k\sum_{\ell=k+1}^n \varepsilon_{\ell}.
\end{equation}
Using the same letter $d$, we next define $d\colon F\rightarrow \mathbb{Z}$ by
$$d(x):=\sum_{i=1}^nd(p_i(x)).$$
Then the equality
\begin{equation}
\label{eq:3-2}
d(xx^{\prime})=d(x)+d(x^{\prime})+[x]\cdot[x^{\prime}]
\end{equation}
holds for $x,x^{\prime}\in F$. Here $[x]\cdot [x^{\prime}]$ is the intersection number
of the homology classes $[x]$ and $[x^{\prime}]$.
Define the map $\tilde{f}\colon \mathcal{N}\times F \rightarrow \mathbb{Z}$ by
\begin{equation}
\label{eq:3-3}
\tilde{f}(\phi,x)=d(\phi(x))-d(x).
\end{equation}
Since the action of $\mathcal{N}$ on $F^{ab}\cong H$ preserves the intersection form,
we see that for each $\phi\in F$,
$\tilde{f}(\phi,\cdot)$ is a homomorphism from $F$ to $\mathbb{Z}$.
Thus the map
$$\tilde{f}\colon \mathcal{N} \rightarrow {\rm Hom}(H,\mathbb{Z})\cong H,$$
which will be also denoted by $\tilde{f}$, is induced. Here the isomorphism
${\rm Hom}(H,\mathbb{Z})\cong H$ is
Poincar\'e duality; $a\in H$ corresponds to the element of ${\rm Hom}(H,\mathbb{Z})$
given by $y\mapsto a\cdot y$.
For simplicity, we write $\rho$ instead of
$\rho \circ \pi \colon \mathcal{N}\rightarrow Sp(H)\cong Sp(2g;\mathbb{Z})$.

\begin{prop}[Morita \cite{M}]
\begin{enumerate}
\item The map $\tilde{f}$ is a twisted 1-cocycle:
$$\tilde{f}(\phi_1\phi_2)=\rho(\phi_2)^{-1}\tilde{f}(\phi_1)+\tilde{f}(\phi_2),\ \phi_1,\phi_2
\in \mathcal{N}.$$
\item Let $f$ be the restriction of $\tilde{f}$ to $\mathcal{M}_{g,1}\subset \mathcal{N}$.
Then in view of the central extension (\ref{eq:1-1}),
$f$ factors through a twisted 1-cocycle $f\colon \mathcal{M}_{g,*}\rightarrow H$ (we use the same letter).
Further, when restricted to the subgroup $\pi_1(\Sigma_g,*)\subset \mathcal{M}_{g,*}$, $f$ coincides
with $(2-2g)$ times the abelianization: $f|_{\pi_1(\Sigma_g,*)}=(2-2g)\theta$.
\end{enumerate}
\end{prop}
In \cite{M}, this cocycle was considered only on $\mathcal{M}_{g,1}$.
But as in the above, it can be naturally defined on $\mathcal{N}$.
In this paper we call $f\colon \mathcal{M}_{g,*}\rightarrow H$ \textit{Morita's twisted 1-cocycle}.

Regard the group $F$ as a subgroup of $\mathcal{N}$ via inner automorphism;
For $x\in F$, the map $F\rightarrow F,\ y\mapsto xyx^{-1}$, is the corresponding element of $\mathcal{N}$.

\begin{lem}
\label{lem:3-2}
For $x\in F$, $\tilde{f}(x)=2[x]$.
\end{lem}

\begin{proof}
Let $x,y\in F$. Using (\ref{eq:3-2}), we compute
\begin{eqnarray*}
\tilde{f}(x,y)&=& d(xyx^{-1})-d(y) \\
&=& d(x)+d(yx^{-1})+[x]\cdot([y]-[x])-d(y) \\
&=& d(x)+d(y)+d(x^{-1})+[y]\cdot[x^{-1}]+[x]\cdot([y]-[x])-d(y) \\
&=& d(x)+d(y)-d(x)+[x]\cdot[x]+[y]\cdot[x^{-1}]+[x]\cdot([y]-[x])-d(y) \\
&=& 2[x]\cdot[y].
\end{eqnarray*}
This proves the lemma.
\end{proof}

\section{Jablow's involution}
We need the value of $\psi$ on a particular element of $\mathcal{M}_{g,*}$.
Consider the endomorphism $\iota$ of $F$ defined by the following:
\begin{equation}
\iota: \begin{cases} A_k \mapsto \left(\prod_{\ell=k}^g [B_g\cdots B_{\ell}A_{\ell},B_{\ell}]B_{\ell}\right)
A_k^{-1}\left(\prod_{\ell=k}^g B_{\ell}^{-1}\right) \\
B_k \mapsto [B_g\cdots B_kA_k,B_k^{-1}]B_k^{-1}
\end{cases}, \ 1\le k\le g.
\label{eq:4-1}
\end{equation}

The following can be checked by a direct computation, so we omit the proof.

\begin{lem}
\label{lem:4-1}
\begin{enumerate}
\item The endomorphism $\iota$ is an involution of $Aut(F)$: $\iota^2={\rm id}$.
\item We have $\iota(\zeta)=(B_g\cdots B_1)\zeta(B_g\cdots B_1)^{-1}$. In particular, $\iota \in \mathcal{N}$.
\end{enumerate}
\end{lem}

The expression of $\iota$ is given in E. Jablow's paper \cite{J} p.231.
Thus we call it \textit{Jablow's involution}.
We write $\pi(\iota)=\bar{\iota}\in \mathcal{M}_{g,*}$. Note that $\rho(\bar{\iota})=-I$.
Topologically, $\bar{\iota}$ is a hyperelliptic involution of $(\Sigma_g,*)$.

\section{Proof of Theorem \ref{main.thm}}
The following two lemmas are needed to prove Theorem \ref{main.thm}: 

\begin{lem}
\label{lem:5-1}
$$\tilde{f}(\iota)=\ ^t(\underbrace{-2,\ldots,-2}_g,-(2g+2),-2g,\ldots,-4).$$
Here we identify $H$ with $\mathbb{Z}^{2g}$ by
$H=\mathbb{Z}\cdot \theta(A_1)\oplus \cdots \oplus \mathbb{Z}\cdot \theta(B_g) \cong \mathbb{Z}^{2g}$.
\end{lem}

\begin{lem}
\label{lem:5-2}
$$\psi(\bar{\iota})=\frac{1}{g-1}\ ^t(\underbrace{1,\ldots,1}_g,-1,-2,\ldots,-g).$$
\end{lem}

The proofs of these lemmas need fairly explicit computations and will be postponed
to the next section. Using these lemmas, we can prove Theorem \ref{main.thm}.

\begin{proof}[Proof of Theorem \ref{main.thm}]
Consider the twisted 1-cocycle $\psi+\frac{1}{2g-2}f\colon \mathcal{M}_{g,*}\rightarrow \frac{1}{2g-2}H$.
The restriction of this cocycle to $\pi_1(\Sigma_g,*)$ is zero. Therefore, by Proposition \ref{prop:1-1},
there exists an element $a\in \frac{1}{2g-2}H$ such that the equation
$$\psi+\frac{1}{2g-2}f=\delta a$$
holds. To determine $a$ we want to evaluate both sides on some element of $\mathcal{M}_{g,*}$.
Unfortunately, $\iota$ is not an element of $\mathcal{M}_{g,1}$ by Lemma \ref{lem:4-1}. So
$f(\bar{\iota})$ may not coincide with $\tilde{f}(\iota)$.
We need to fix $\iota$ to obtain an element of $\mathcal{M}_{g,1}$.
Let $x_B:=(B_g\cdots B_1)^{-1}\in F\subset \mathcal{N}$. By the second statement of Lemma \ref{lem:4-1},
$x_B\cdot \iota \in \mathcal{M}_{g,1}$.

Now we compute
\begin{eqnarray*}
f(\pi(x_B\cdot \iota)) &=& \tilde{f}(x_B\cdot \iota) \\
&=& \rho(\iota)^{-1}\tilde{f}(x_B)+\tilde{f}(\iota) \\
&=& -\ ^t(\underbrace{0,\ldots,0}_g,\underbrace{-2,\ldots,-2}_g)
+\ ^t(\underbrace{-2,\ldots,-2}_g,-(2g+2),-2g,\ldots,-4) \\
&=& \ ^t(\underbrace{-2,\ldots,-2}_g,-2g,-2g+2,\ldots,-2)
\end{eqnarray*}
by using Lemma \ref{lem:3-2} and Lemma \ref{lem:5-1}, and
\begin{eqnarray*}
\psi(\pi(x_B\cdot \iota))&=&\psi(\pi(x_B)\cdot \bar{\iota}) \\
&=& \rho(\bar{\iota})^{-1}\psi(\pi(x_B))+\psi(\bar{\iota}) \\
&=& -\ ^t(\underbrace{0,\ldots,0}_g,\underbrace{-1,\ldots,-1}_g)
+\frac{1}{g-1}\ ^t(\underbrace{1,\ldots,1}_g,-1,-2,\ldots,-g) \\
&=& \frac{1}{g-1}\ ^t(\underbrace{1,\ldots,1}_g,g-2,g-3,\ldots,-1)
\end{eqnarray*}
by using Lemma \ref{lem:5-2} and Proposition \ref{prop:2-1}.
On the otherhand, we compute
$$\left(\psi+\frac{1}{2g-2}f\right) (\pi(x_B\cdot \iota))=
\delta a(\pi(x_B\cdot \iota))=\rho(x_B\cdot \iota)^{-1}\cdot a-a=-2a.$$
Combining these computations together, we obtain $a=a_0$.
\end{proof}

\section{Proofs of Lemma \ref{lem:5-1} and Lemma \ref{lem:5-2}}
In this section we prove Lemma \ref{lem:5-1} and Lemma \ref{lem:5-2}.

\begin{proof}[Proof of Lemma \ref{lem:5-1}]
Let $k\in \{1,\ldots,g\}$.
By (\ref{eq:4-1}), we have
$$p_i(\iota(A_k))=
\begin{cases}
1 & (i<k) \\
\beta \alpha \beta \alpha^{-1} \beta^{-1} \alpha^{-1} \beta^{-1} & (i=k) \\
\beta \alpha \beta \alpha^{-1} \beta^{-1} \beta^{-1} & (i>k)
\end{cases}$$
and
$$p_i(\iota(B_k))=
\begin{cases}
1 & (i\neq k) \\
\beta \alpha \beta^{-1} \alpha^{-1} \beta^{-1} & (i=k).
\end{cases}$$
Direct computations according to (\ref{eq:3-1}) shows that
$$d(\beta \alpha \beta \alpha^{-1} \beta^{-1} \alpha^{-1} \beta^{-1})=4,\ 
d(\beta \alpha \beta \alpha^{-1} \beta^{-1} \beta^{-1})=2,\ 
d(\beta \alpha \beta^{-1} \alpha^{-1} \beta^{-1})=-2.$$
From these, we have
$\tilde{f}(\iota,A_k)=4+2(g-k)$ and $\tilde{f}(\iota,B_k)=-2$. This proves the lemma.
\end{proof}

We next proceed to the proof of Lemma \ref{lem:5-2}.
At first glance it seems difficult to compute
the vector $w=w(\phi,[X,p,f])$ from the definition (\ref{eq:2-1})
because it is not easy to compute quadratic period classes in general.
However, if $\phi$ is a hyperelliptic involution such as $\bar{\iota}$,
all terms concerning quadratic period classes cancel and disappear
as we will see in the below for the case $\phi=\bar{\iota}$. 

\begin{proof}[Proof of Lemma \ref{lem:5-2}]
Let us introduce a complex structure on $(\Sigma_g,*)$ and denote it by $(X_0,p_0)$.
Consider the point $[X_0,p_0,{\rm id}]\in \mathcal{T}_{g,1}$.
According to (\ref{eq:2-1}),
$$w=w(\bar{\iota},[X_0,p_0,{\rm id}])=-\eta([X_0,p_0,\bar{\iota}])-\eta([X_0,p_0,{\rm id}])$$
since ${\bar{\iota}}^{-1}=\bar{\iota}$ and $A(\bar{\iota},[X_0,p_0,{\rm id}])=-I$.
Thus the $j$-th component of the vector $w$ is given by
\begin{equation}
\label{eq:6-1}
w_j=\frac{1}{g-1}\left( -\tau_{jj}+\sum_{k=1}^g (Q_{jk}(A_k)+Q_{jk}(\bar{\iota}_*A_k)) \right).
\end{equation}
Now quadratic period classes satisfy the following property (see for instance, \cite{J} p.222):
$$Q_{jk}(T_1T_2)=Q_{jk}(T_1)+Q_{jk}(T_2)+\omega_j(T_1)\omega_k(T_2),\ T_1,T_2\in \pi_1(X_0,p_0),$$
where $\omega(T)=\int_{[T]}\omega$ is the integration of the holomorphic 1-form $\omega$ along
the homology class $[T]$. From this, the following formulae are easily derived;
$$Q_{jk}(T_1\cdots T_m)=\sum_{n=1}^m Q_{jk}(T_n)+
\sum_{1\le n < n^{\prime} \le m}\omega_j(T_n)\omega_k(T_{n^{\prime}});$$
$$Q_{jk}([S,T])=\omega_j(S)\omega_k(T)-\omega_j(T)\omega_k(S);$$
$$Q_{jk}(T^{-1})=\omega_j(T)\omega_k(T)-Q_{jk}(T),$$
where $S,T,T_i\in \pi_1(X_0,p_0)$.
We will use these formulae freely in the following computations.
Set $E_{\ell}:=[B_g\cdots B_{\ell}A_{\ell},B_{\ell}]B_{\ell}$. Then
$\bar{\iota}_*A_k=E_kE_{k+1}\cdots E_gA_k^{-1}B_k^{-1}\cdots B_g^{-1}$.
We compute
\begin{eqnarray*}
Q_{jk}(\bar{\iota}_*A_k)&=& Q_{jk}(E_kE_{k+1}\cdots E_gA_k^{-1}B_k^{-1}\cdots B_g^{-1}) \\
&=& \sum_{\ell=k}^g Q_{jk}(E_{\ell})+Q_{jk}(A_k^{-1})+\sum_{\ell=k}^g Q_{jk}(B_{\ell}^{-1})
+\sum_{k\le \ell < \ell^{\prime}} \omega_j(E_{\ell})\omega_k(E_{\ell^{\prime}}) \\
&\ &+\sum_{\ell=k}^g \omega_j(E_{\ell})\omega_k(A_k^{-1})
+\sum_{k\le \ell,k\le \ell^{\prime}}\omega_j(E_{\ell})\omega_k(B_{\ell^{\prime}}^{-1})
+\sum_{\ell=k}^g \omega_j(A_k^{-1})\omega_k(B_{\ell}^{-1}) \\
&\ &+\sum_{k \le \ell < \ell^{\prime}} \omega_j(B_{\ell}^{-1})\omega_k(B_{\ell^{\prime}}^{-1});
\end{eqnarray*}

\begin{eqnarray*}
Q_{jk}(E_{\ell})&=&Q_{jk}([B_g\cdots B_{\ell}A_{\ell},B_{\ell}]B_{\ell}) \\
&=& Q_{jk}([B_g\cdots B_{\ell}A_{\ell},B_{\ell}])+Q_{jk}(B_{\ell}) \\
&=& \omega_j(B_g\cdots B_{\ell}A_{\ell})\omega_k(B_{\ell})-\omega_j(B_{\ell})\omega_k(B_g\cdots B_{\ell}A_{\ell})
+Q_{jk}(B_{\ell}) \\
&=& \left(\sum_{m=\ell}^g \tau_{jm}+\delta_{j\ell}\right)\tau_{k\ell}
-\tau_{j\ell}\left(\sum_{m=\ell}^g \tau_{km}+\delta_{k\ell}\right)
+Q_{jk}(B_{\ell});
\end{eqnarray*}

$$Q_{jk}(A_k^{-1})=\omega_j(A_k)\omega_k(A_k)-Q_{jk}(A_k)=\delta_{jk}-Q_{jk}(A_k);$$
$$Q_{jk}(B_{\ell}^{-1})=\omega_j(B_{\ell})\omega_k(B_{\ell})-Q_{jk}(B_{\ell})
=\tau_{j\ell}\tau_{k\ell}-Q_{jk}(B_{\ell});$$
$\omega_j(E_{\ell})=\tau_{j\ell}$, and $\omega_k(E_{\ell})=\tau_{k\ell}$. Therefore,
$Q_{jk}(A_k)+Q_{jk}(\bar{\iota}_*A_k)$ is equal to
\begin{eqnarray}
&\ & \sum_{\ell=k}^g \left( \left(\sum_{m=\ell}^g \tau_{jm}+\delta_{j\ell}\right) \tau_{k\ell}
-\tau_{j\ell}\left( \sum_{m=\ell}^g \tau_{km}+\delta_{k\ell}\right) \right)
+\delta_{jk}+\sum_{\ell=k}^g \tau_{j\ell}\tau_{k\ell} \nonumber \\
&\ & +\sum_{k\le \ell < \ell^{\prime}}\tau_{j\ell}\tau_{k\ell^{\prime}}-\sum_{\ell=k}^g \tau_{j\ell}
-\sum_{k\le \ell,k\le \ell^{\prime}}\tau_{j\ell}\tau_{k\ell^{\prime}}
+\sum_{\ell=k}^g \delta_{jk}\tau_{k\ell}
+\sum_{k\le \ell < \ell^{\prime}} \tau_{j\ell}\tau_{k\ell^{\prime}}. \label{eq:6-2}
\end{eqnarray}
Now no quadratic period classes appear here. We observe that the degree two terms with respect
to $\tau$ (e.g. $\tau_{jm}\tau_{k\ell}$) appearing in (\ref{eq:6-2}) all cancel out.
Thus (\ref{eq:6-2}) is equal to
$$ \sum_{k\le \ell}(\delta_{j\ell}\tau_{k\ell}-\delta_{k\ell}\tau_{j\ell})+\delta_{jk}
-\sum_{\ell=k}^g \tau_{j\ell}+\sum_{\ell=k}^g \delta_{jk}\tau_{k\ell},$$
and we compute $\sum_{k=1}^g \left( Q_{jk}(A_k)+Q_{jk}(\bar{\iota}_*A_k) \right)$ as
\begin{eqnarray}
&\ &\sum_{1\le k\le \ell \le g}(\delta_{j\ell}\tau_{k\ell}-\delta_{k\ell}\tau_{j\ell})
+\sum_{k=1}^g \delta_{jk}-\sum_{1\le k\le \ell \le g}\tau_{j\ell}
+\sum_{1\le k\le \ell \le g} \delta_{jk}\tau_{k\ell} \nonumber \\
&=& \sum_{1\le k\le j}\tau_{kj}-\sum_{k=1}^g\tau_{jk}+1
-\sum_{\ell=1}^g \ell \tau_{j\ell}+\sum_{j\le \ell \le g}\tau_{j\ell} \nonumber \\
&=& 1+\tau_{jj}-\sum_{\ell=1}^g \ell \tau_{j\ell}. \label{eq:6-3}
\end{eqnarray}
Here we use the symmetry condition $\tau_{kj}=\tau_{jk}$ of period matrices.
Substituting (\ref{eq:6-3}) into (\ref{eq:6-1}) we obtain
$$w_j=\frac{1}{g-1}\left( 1-\sum_{\ell=1}^g \ell \tau_{j\ell}\right),$$
and
$$w=(I,\tau)\cdot \frac{1}{g-1}\ ^t(\underbrace{1,\ldots,1}_g,-1,-2,\ldots,-g).$$
Since the $2g$ vectors of the colums of the $g\times 2g$ matrix $(I,\tau)$ are
linearly independent over the real numbers,
the lemma follows.
\end{proof}

\noindent \textbf{Acknowledgments:}
I am grateful to Professor Nariya Kawazumi for telling me about Earle's paper,
and Professor Shigeyuki Morita for useful comments about twisted 1-cocycles on
the mapping class groups. I also thank them for their warm encouragements.  
This research is supported by JSPS Research
Fellowships for Young Scientists (19$\cdot$5472).

\noindent \textsc{Yusuke Kuno\\
Graduate School of Mathematical Sciences,\\
The University of Tokyo,\\
3-8-1 Komaba Meguro-ku Tokyo 153-0041, JAPAN}

\noindent \texttt{E-mail address:kunotti@ms.u-tokyo.ac.jp}

\end{document}